# A model of discrete Kolmogorov-type competitive interaction in a two-species ecosystem

Sudeepto Bhattacharya[1], L.M. Saha[2]

([1]Department of Mathematics, [2]IIIMIT, School of Natural Sciences, Shiv Nadar University,

P.O. Shiv Nadar University, Greater Noida, Gautam Buddha Nagar, Uttar Pradesh, India)

An ecosystem is a nonlinear dynamical system, its orbits giving rise to the observed complexity in the system. The diverse components of the ecosystem interact in discrete time to give rise to emergent features that determine the trajectory of system's time evolution.

The paper studies the evolutionary dynamics of a toy two-species ecosystem modelled as a discrete-time Kolmogorov system. It is assumed that only the two species comprise the ecosystem and compete with each other for obtaining growth resources, mediated through inter- as well as intra-specific coupling constants to obtain resources for growth. Numerical simulations reveal the transition from regular to irregular dynamics and emergence of chaos during the process of evolution of these populations. We find that the presence or absence of chaotic dynamics is being determined by the interspecific interaction coefficients. For values of the interspecific interaction constants widely apart, the system emerges to regular dynamics, implying coexistence of the competing populations on a long- term evolutionary scenario.

**AMSC:** 92D40

## I. Introduction

An ecosystem is a discrete map defined on a compact, normed space. The various states of ecosystem that are observed in nature are the orbits of the map. Following common usage of the term, we shall henceforth refer to an orbit of the above map as an ecosystem.

A given ecosystem evolves in time by evolving complex structures and patterns due to the emergent behaviour of the interactions between the components comprising it, qualifying as a complex adaptive system. These patterns, in turn, serve as indicators of long-term coexistence or extinction of the constituent species populations. This feature of the ecosystem qualifies it to be a complex adaptive system [1, 2, 3, 4, 5].

The question of dynamic coexistence of competing species in a complex ecosystem has been a problem of prime importance in dynamical systems research, attracting a significant body of studies [6, 7, 8, 9, 10, 11, 12]. Both, intra-and inter-specific competitions constitute one of the fundamental classes of ecological interaction, leading to possible coexistence or extinction of competing species in a given ecosystem. Hence, competitive interactions are important objects of study in order to gain an understanding of the dynamics of evolution and emergence of complex structures in the ecosystem.

The objective of the present work is to explore the discrete-time dynamical behaviour of a competitive two-species "toy" model of an ecosystem such that both the species occupy the same trophic level. We use numerical experiments with an aim to describe and understand some possible characteristic features of the emergence of dynamical regimes and hence of complexity in the evolving ecosystem. The work is a descriptive modelling of the "toy" ecosystem behaviour. In particular, we study the explicit influence of the competitive interaction on a long-term coexistence of the species, which have non-overlapping generations. We assume that the "toy" ecosystem comprises populations of only two species at the same trophic level, which we call species A and B respectively, such that the length of their generations is similar. The paper makes a further simplification by assuming that apart from these two species, populations of no other species, floral or faunal, is present in the ecosystem, meaning thereby that the ecosystem entirely comprises species A and B only.

It may be noted that discrete-time complex dynamics of ecosystems have been well-studied since long. Such systems have been analysed by authors, primarily with a focus on the predator-prey communities or on a single species population [7, 8, 9]. The work presented here has its focus instead on the inter-specific competition that may best be conceptualized at the same trophic level of the ecosystem.

An appropriate model for addressing our objective is the discrete-time Kolmogorov system, which adequately describes the population dynamics two-species population with non-overlapping generations [10, 11, 12, 13, 14, 15, 16].

In this work we study the time evolution of non-overlapping generations of the two species A and B, comprising a typical "toy" ecosystem. Without loss of generality, we choose to study the competitive response of population of A to population of B using numerical simulations to understand the various discrete dynamical signatures of such an interaction, and address the question of coexistence of these populations. In section II, we propose the mathematical model of discrete-time Kolmogorv system governing the dynamics of the toy ecosystem and pose the research question. In section III, we report sample results of the numerical experiments and simulations performed on the proposed model, in order to address and answer the question posed, and conclude the paper in section IV, after drawing inferences from the results obtained.

## II. Modelling

Let $(X, \|\cdot\|)$ be a compact, normed space. Let $\bar{x}_0 \in X$. Let

$$T : X \to X \tag{1}$$

be a discrete-time map, and let $\bar{x}_n := T^n(x_0)$ be the iterates of $T$ for $n \in N = \{0,1,...\}$, obtained as orbits on discrete time evolution of $T$. The map $T$ is the ecosystem in this work, with iterates denoting different states of the ecosystem in the state space $X$.

Let the state space of the toy ecosystem be the normed space $(X, \|\cdot\|)$. Assume that the population dynamics of the two competing species are governed by the discrete-time Kolmogorov map having the form

$$x_{t+1}^{(i)} = x_t^{(i)} T^{(i)}(\bar{x}_t), i = 1,2, t \in N, \tag{2}$$

In the system (2) above, we consider $T$ to be a two-dimensional logistic map and denote the population densities of the B and A species at time (generation) $t = n$ respectively by $x_n^{(i)}$, $i = 1,2$. The corresponding vector of population densities is given by $\bar{x}_n = \begin{bmatrix} x_n^{(1)} \\ x_n^{(2)} \end{bmatrix} \in X$, $x_0^{(i)} > 0, i = 1,2$.

In our study, the state space of the system (2) is a subspace of $\Re^2$, the two-dimensional real vector space.

For the purpose of our modelling, we write the system (2) explicitly as

$$\left. \begin{array}{l} x_{n+1}^{(1)} = x_n^{(1)}[r_1(1 - c_1 x_n^{(1)} - c_2 x_n^{(2)})] \\ x_{n+1}^{(2)} = x_n^{(2)}[r_2(1 - c_3 x_n^{(1)} - c_4 x_n^{(2)})] \end{array} \right\}, \tag{3}$$

where $r_1 \geq 0, r_2 \geq 0$ are respectively the logistic growth rates of the populations of A and B species, and the non-negative coupling constants $c_1, c_4$ are the intra-species competition coefficients while $c_2, c_3$ are the corresponding quantities representing inter-species competition [17].

For the sake of simplicity of notation, we shall write $x_n^{(1)}$ as $x_n$ and $x_n^{(2)}$ as $y_n$. Since our concern in the present work is with the A species population's response to competitive interaction by the B population, we shall assume a simplification and study the discrete-time Kolmogorov dynamical system (3) with only $r_2$ as the control parameter in the numerical simulations. With such a simplifying assumption, the research question for our work may be stated as: *what are the asymptotic behaviours of the trajectories of the Kolmogorov system (3), when the growth dynamics is controlled by the parameter $r_2$?*

We address the question posed above by performing numerical experiments and corresponding simulations of the dynamics of the ecosystem (3) and report sample instances of the results of the simulation in this paper.

### III. Numerical experiments, simulations and results

In the following, we report the experimental results obtained by numerically simulating the discrete Kolmogorov system (3), to investigate (I) stability of the system; (II) bifurcations in the system and (III) Lyapunov exponents of the system.

### (I) Stability of the system

To examine the stability of the system (3) at a given fixed point $\bar{x}^* = \begin{bmatrix} x^* \\ y^* \end{bmatrix} = (x^*, y^*)$, we evaluate the Jacobian matrix $J$ of the map $T$ at this point and then subject $J$ to standard test procedure for obtaining the conditions for existence of (stable or unstable) fixed points of the map.

We begin by setting the coupling constants $c_2 = c_3 = 0$ for the first part of this investigation. The assumption makes the coupled system (3) get decoupled, and the two constituents behave as two independent equations, and hence the individual species population densities do not influence one another. The following table summarizes the results obtained:

| S. No | Fixed point | Parameter values required for stability |
|---|---|---|
| 1 | $(0,0)$ | $0 < r_1 < 1$, $0 < r_2 < 1$, $c_1, c_4$ arbitrary |
| 2 | $\left(\dfrac{r_1 - 1}{c_1 r_1}, 0\right)$ | $|r_2| < 1$, $|r_1 - 2| < c_1 r_1$ |
| 3 | $\left(0, \dfrac{r_2 - 1}{c_4 r_2}\right)$ | $|r_1| < 1$, $|r_2 - 2| < c_4 r_2$ |
| 4 | $\left(\dfrac{r_1 - 1}{c_1 r_1}, \dfrac{r_2 - 1}{c_4 r_2}\right)$ | $|r_1 - 1| < c_1 r_1$, $|r_2 - 1| < c_4 r_2$ |

We then carry out the fixed point analysis in presence of the inter-species interactions. For this part of the exercise we assign non-zero arbitrary but feasible values to the coupling constants, including $c_2$ and $c_3$. Due to the reason mentioned in the foregoing, our study in this work is concerned about the dependence of dynamics of the system (3) on the parameter $r_2$, the growth rate of the A ecosystem. We therefore use $r_2$ as the control parameter in the system dynamics and treat all other parameters as constants for any given instance. It may be noted that because of being logistic in character, $r_2 \in [0,4]$. By changing the parameter values suitably, coordinates of the fixed points for the map (2) could be obtained. The above stability analysis also indicated the feasible ranges for the permissible values of the various parameters involved in the model.

*(II) Bifurcations in the system*

Keeping all other parameters of the system (3) at constant values for each instance, we vary $r_2$ to obtain the bifurcation diagrams of the system dynamics given by the map $T$. We observe from the diagrams given in Fig.1 that the dynamics of the system exhibits regularity for

$$r_1 = 3.0, c_1 = 1.8, c_2 = \{0.1, 0.2, 0.3\}, c_3 = 0.6, c_4 = 2.5, r_2 \in [2.8, 4.0]$$

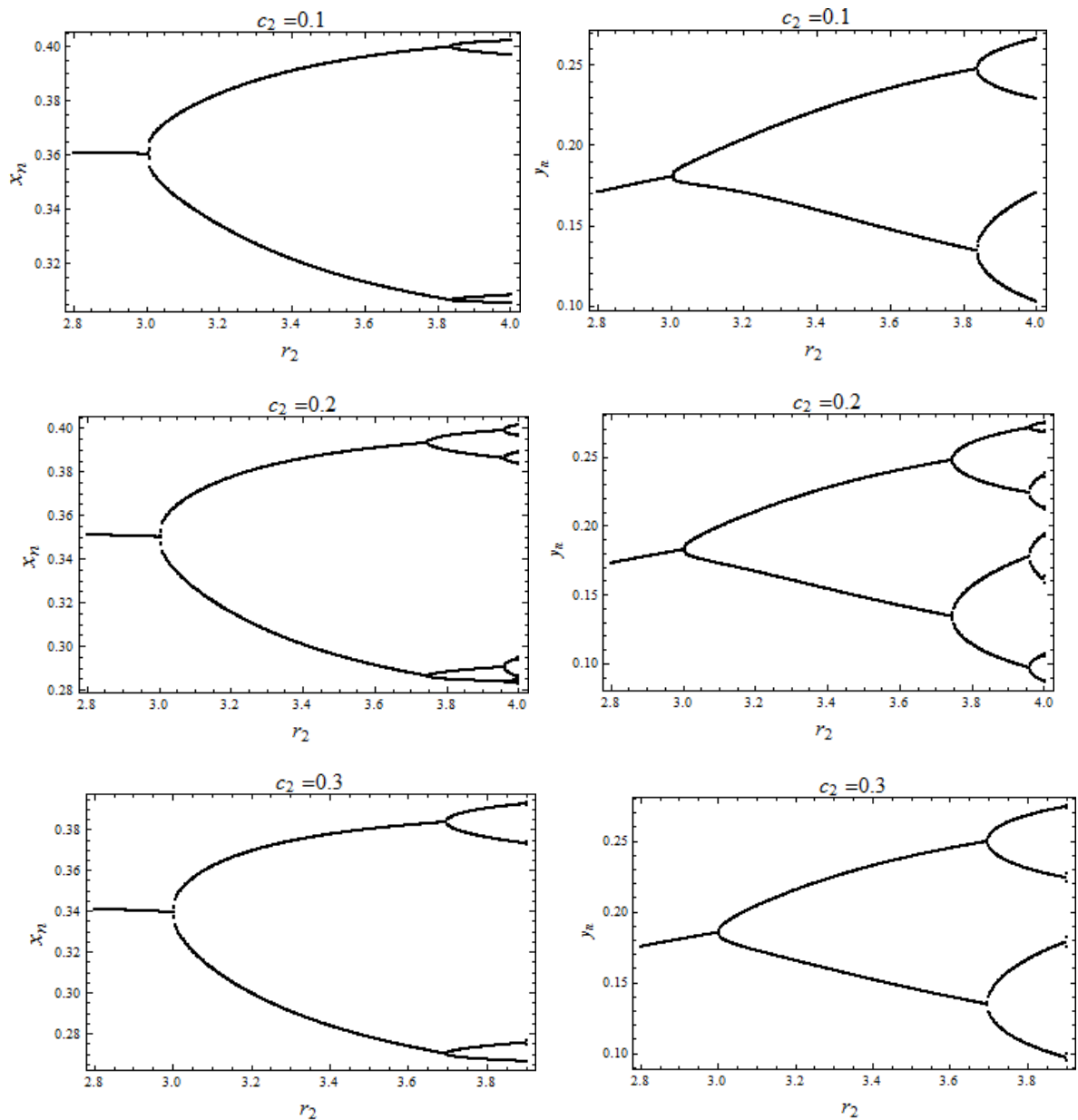

**Fig. 1 Regular behaviour of the system dynamics with variation of $r_2$ for different values of coupling coefficients $c_2$, when other parameters are $r_1$ = 3.0, $c_1$ = 1.8, $c_3$ = 0.6, $c_4$ = 2.5. Initial conditions are taken as ($x_0^{(1)}$, $x_0^{(2)}$) = (0.2, 0.1), screenshot after 400 iterations.**

We however observe from the simulation that when the values of $c_2$ and $c_3$ are chosen closer to one another, the dynamics of the system tends to behave chaotically as could be seen from Fig. 2 below:

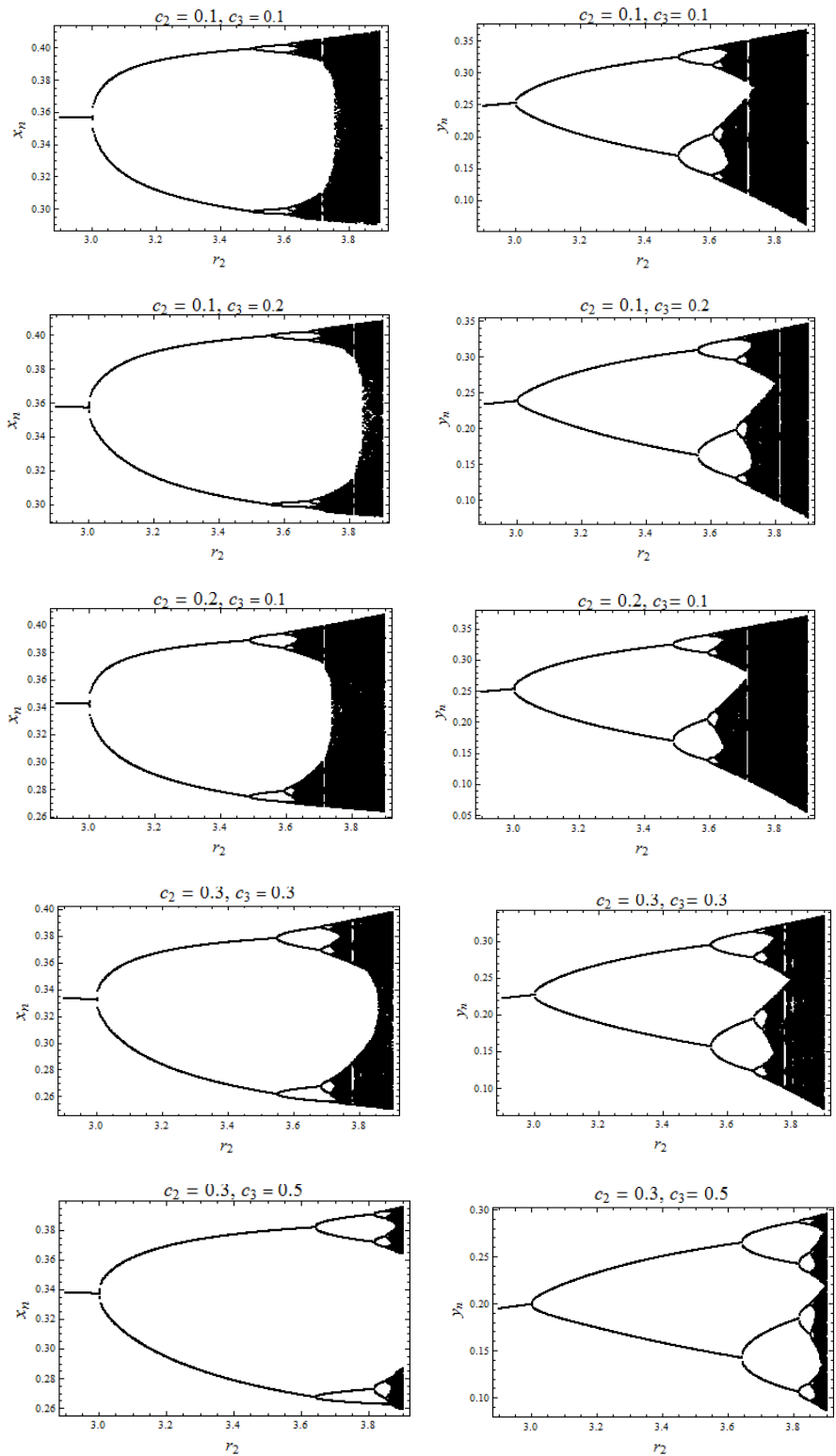

**Fig. 2 Chaotic behaviour of the system dynamics with variation of $r_2$ for closer values of coupling coefficients $c_2, c_3$ when other parameters are $r_1 = 3.0$, $c_1 = 1.8$, $c_4 = 2.5$ and the initial conditions are taken as ($x_0^{(1)}$, $x_0^{(2)}$) = (0.2, 0.1). Screenshot after 400 iterations.**

It may thus be inferred from the foregoing that proximity in the values of inter-specific interaction coefficients serves as one of the determinants of emergence of chaos in the population dynamics of the species.

*(III) Lyapunov exponents for the system*

Lyapunov exponent is obtained as an average measure of the exponential divergence of two orbits evolving under the given map $T$, initiated with infinitesimal separation [18, 20]. It therefore provides a measure of local (between neighbouring points) average expansions of the system (3) as the system evolves in time starting from the initial state at time $t = 0$. Among various characteristics of chaotic dynamics, positivity of the largest Lyapunov exponent happens to be a principal one. We shall call the orbit set of $T$ as chaotic if its asymptotic measure has positive Lyapunov exponent [18, 19, 20, 21, 22].

Let $U_i = \{\bar{x}_{i-1}, \bar{x}_{i+1}\}$, $i = 0,...,n-1$ be the set of all neighbours of the point $\bar{x} \in X$. Then, we may write $U_0 = \|\bar{x}_0 - \bar{y}_0\|$ for $\bar{x}_0, \bar{y}_0 \in X$. Lyapunov exponent for the map $T$ is then given by

$$\lambda_T = \lim_{n \to \infty} \frac{1}{n} \ln \left\| \prod_{i=0}^{n-1} J(\bar{x}_i) U_0 \right\| \tag{4}$$

and, $\|\bar{x}_n - \bar{y}_n\| \approx e^{\lambda_T n}$ describes the exponential separation with iteration $n$ between two neighbouring points that evolve under the map $T$, with infinitesimal initial separation.

The following figures, sampled from the results obtained through numerical simulations, depict the phase portraits of the discrete time evolution of the two competing species, alongside the functional dependence of the Lyapunov exponent $\lambda\, (= \lambda_T)$ on the number of iterations $n$. It may be observed from the figures that the behaviour of the system is in consonance with the results and of *(B)* illustrated in figs 1 and 2 above, and the inference drawn.

The diagrams in Fig 3 show that for values of the inter-specific interaction coefficients that are chosen apart from one another, the two species exhibit long-term co-existence in the toy ecosystem. The Lyapunov exponents indicate that the system dynamics enters chaotic regime when the above coefficients are chosen closer or identical in value to one another.

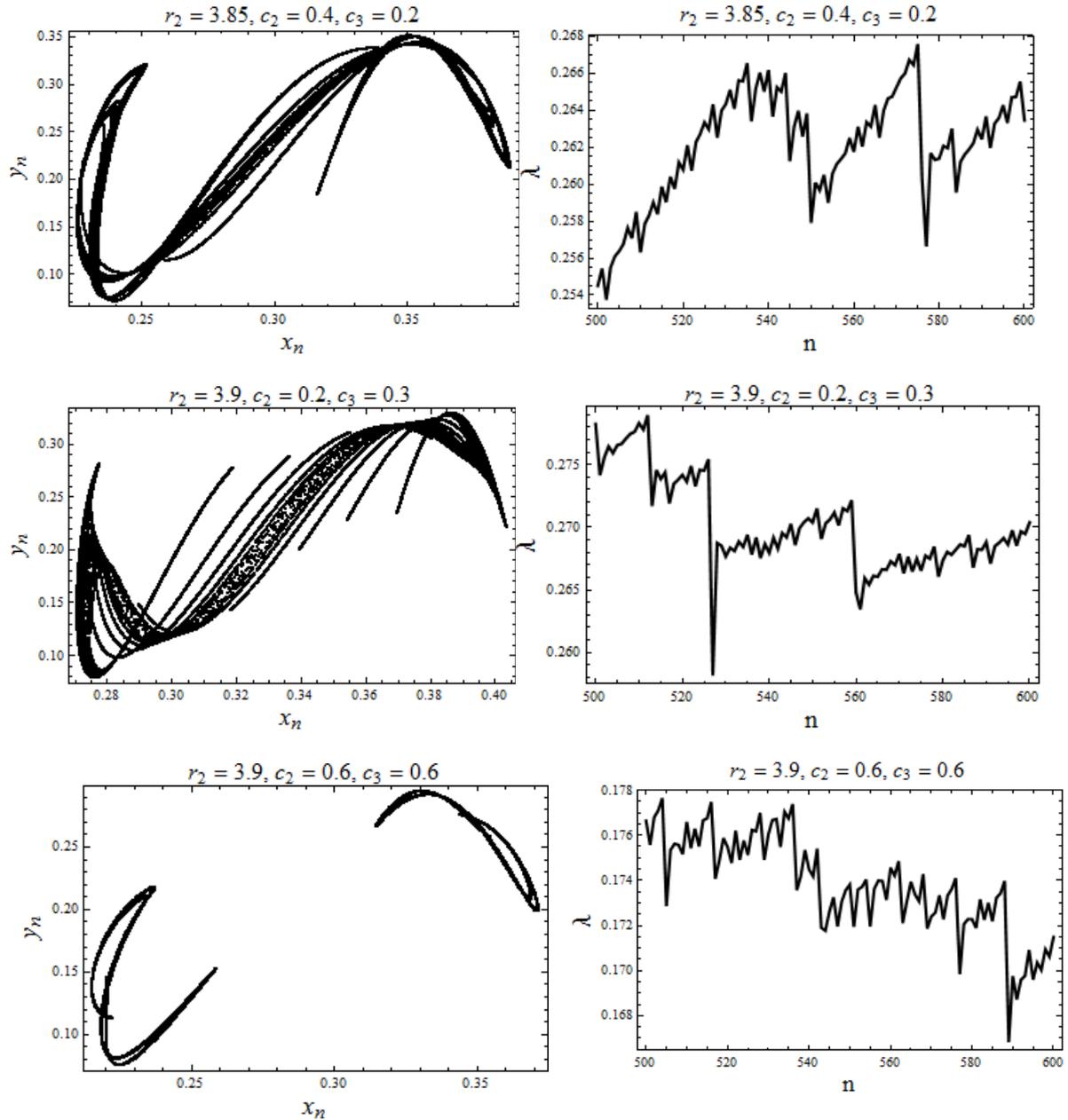

**Fig.3 Phase diagrams and corresponding Lyapunov exponents for the system dynamics. Other parameters are $r_1 = 3.0$, $c_1 = 1.8$, $c_5 = 2.5$ and initial conditions are $(x_0^{(1)}, x_0^{(2)}) = (0.1, 0.1)$. All plots are between iteration numbers 500 and 600.**

## IV Conclusion

Three related aspects of the dynamics of the complex, discrete-time Kolmogorov system (3) were studied in this paper through numerical simulation. We find the inter-species coupling coefficients $c_2$ and $c_3$ to be important determinants in the long-term dynamics and consequent emergence of patterns in the "toy" model. When their values are chosen in close proximity to each other, or are chosen to be identical, the dynamics of the system (3) enters a chaotic regime, indicated both, by the bifurcation diagrams and the Lyapunov exponents for the system. Existence of chaos in the evolution of populations of the competing species may therefore be considered as an emergent feature in the long-term behaviour of the ecosystem.

In this context, it is worth noting that a discrete dynamical model of competitive outcomes mentions emergent dynamical scenarios even in contradiction to the well-known Competitive Exclusion Principle [23].

However, we observe for values of the above coefficients chosen distant apart, the dynamical behaviour of the competitive outcomes remain regular. It may thus be conjectured that the two components of the ecosystem, namely the populations of species A and B, would continue to competitively coexist through the long-term evolution of the system (3), provided that the intra-specific competitive interactions have an appreciable difference in strengths, thus yielding the emergent phenomenon of population coexistence in the long-term dynamics of the ecosystem.

**References**

[1] Anand, M., Gonzales, A.,Guichard, F., Kolasa, J., Parrott, L., 2010. Ecological systems as complex systems: challenges for an emerging science. Diversity **2**: 395 – 410.

[2] Gell-Mann, M., 1994. Complex adaptive systems. Complexity: Metaphors, models and reality Proc., XIX: 17 – 45, Addison-Wesley.

[3] Levin, S. A., 1998a. Ecosystems and biosphere as complex adaptive systems. Ecosystems 1: 431 – 436.

[4] Levin, S. A., 2002. Complex adaptive systems: exploring the known, the unknown and the unknowable. Bull. (New series) of the Amer. Math. Soc. 40 (1): 3 – 19.

[5] Slobodkin, L. B., 1992. Simplicity and Complexity in Games of the Intellect, Harvard University Press, Cambridge, Massachusetts.

[6] Beddington, J. R., Free, C. A., Lawton, J. H., 1975. Dynamic complexity in predator-prey models framed in difference equations. Nature 255: 58 – 60.

[7] Ricker, W. E., 1954. Stock and recruitment. J. Fes. Res. Board Can. 11(5): 559 – 623.

[8] May, R. M., 1974. Stability and Complexity in Model Ecosystems. Princeton University Press, Princeton.

[9] Frisman, E. Y., Neverova, G.P., RevutskayaO. L., 2011. Complex dynamics of the population with a simple age structure. Ecological Modelling 222: 1943 – 1950.

[10] Franke, J. E., Yakubu, A. A., 1991. Mutual exclusion versus coexistence for discrete competitive systems. J. Math. Bio. 30: 161 – 168.

[11] Hassell, M. P., 2000. The Spatial and Temporal Dynamics of Host-Parasitoid Interactions, Oxford University Press, Oxford.

[12] Hofbauer, J., Hutson, V., Jansen, W., 1987. Coexistence for systems governed by difference equations of Lotka - Volterra type. J. Math. Bio. 25: 553 – 570.

[13] Hutson, V., Moran, W., 1982. Persistence of species obeying difference equations. J. Math. Bio. 15: 203 – 213.

[14] Kot, M., 2002. Elements of Mathematical Ecology, Cambridge University Press, Cambridge.

[15] Bhattacharya, S., Saha, L. M., 2013. Dynamics in a two-species Kolmogorov system. International J. Engineering and Innovative Technology 2(8): 92 – 98.

[16] Kon, R., 2004. Permanence of discrete-time Kolmogorov systems for two species and saturated fixed points. J. Math. Bio. 48: 57 – 81.

[17] Ludwig, D., Walker, B., Holling, C. S., 1997. Sustainability, stability and resilience. Conservation Ecology 1(1)7: 1 – 24.

[18] Amigo, J. M., Kocarev, L., Szczepanski, J., 2008. On some properties of the discrete Lyapunov exponent. Phys. Lett. A 372: 6265 – 6268.

[19] Geist, K., Parlitz, U., Lauterborn, W., 1990. Comparison of different methods for computing Lyapunov exponents. Progr. Theor. Phys. 83, 875 – 893.

[20] Hirsch, M. W., Smale, S., Devaney, R. L. 2004. Differential Equations, Dynamical Systems and an Introduction to Chaos, Academic Press, San Diego.

[21] Jost, J., 2005. Dynamical Systems Examples of Complex Behavior, Springer, Berlin, Heidelberg, New York.

[22] Saha, L. M., Prasad, S., Yuasa, M., 2012. Measuring chaos, topological entropy and correlation dimension in discrete maps. Science and Technology 24: 11 – 23, Kinki University, Japan.

[23] Cushing, J. M., LeVarge S., 2005. Some discrete competition models and the principle of competitive exclusion. In Difference equations and discrete dynamical systems. Proceedings of the Ninth International Conference, pp: 283 - 302, Allen L. J. S., Aulbach B., Elyadi S., Sacker R (Eds.) World Scientific, Singapore.